\documentclass[5p,times]{elsarticle}
\usepackage{algorithm}
\usepackage{algpseudocode}
\usepackage[colorinlistoftodos,prependcaption,textsize=tiny]{todonotes}
\usepackage{listings}
\usepackage{booktabs}
\usepackage[position=bottom]{subcaption}

\usepackage{hyperref}
\hypersetup{breaklinks=true}

\definecolor{clr-background}{RGB}{255,255,255}
\definecolor{clr-text}{RGB}{0,0,0}
\definecolor{clr-string}{RGB}{163,21,21}
\definecolor{clr-namespace}{RGB}{0,0,0}
\definecolor{clr-preprocessor}{RGB}{255,128,0}
\definecolor{clr-preprocessor2}{RGB}{153,76,0}
\definecolor{clr-preprocessor3}{RGB}{102,0,204}
\definecolor{clr-keyword}{RGB}{0,0,255}
\definecolor{clr-type}{RGB}{43,145,175}
\definecolor{clr-variable}{RGB}{0,0,0}
\definecolor{clr-constant}{RGB}{111,0,138} % macro color
\definecolor{clr-comment}{RGB}{0,128,0}
\definecolor{light-gray}{gray}{0.95}
\lstdefinestyle{cstyle}{
    language=C++,
    xleftmargin=.05\linewidth,
	backgroundcolor=\color{white},
	basicstyle=\scriptsize\ttfamily\color{clr-text}, % any text
	stringstyle=\color{clr-string},
	commentstyle=\color{clr-comment},
    keywordstyle = {\color{clr-type}},
	tabsize=4,
	captionpos=b,
    columns=fullflexible,
	basewidth={0.6em,0.55em},
	moredelim=[is][\color{clr-preprocessor}]{|>}{<|}	
}
\definecolor{commentsColor}{rgb}{0.0, 0.4, 0.2}
\definecolor{keywordsColor}{rgb}{0.6, 0.0, 0.4}
\definecolor{stringColor}{rgb}{0.0, 0.0, 0.7}

\newcommand{\ignore}[1]{}
\begin{document}

\begin{frontmatter}
\title{ShyLU-node: On-node Scalable Solvers and Preconditioners\\
       Recent Progresses and Current Performance}

\author[snl]{Ichitaro Yamazaki}
%\email{iyamaza@sandia.gov}
\author[snl]{Nathan Ellingwood}
%\email{ndellin@sandia.gov}
\author[snl]{Sivasankaran Rajamanickam}
%\email{srajama@sandia.gov}

\fntext[snl]{Sandia National Laboratories, Albuquerque, New Mexico, U.S.A}

\begin{abstract}

ShyLU-node is an open-source software package that implements linear solvers and preconditioners on shared-memory multicore CPUs or on a GPU. 
These ShyLU solvers and preconditioner can be used as a stand-alone global problem solver, as a local subdomain solver for a domain decomposition (DD) preconditioner, or as a coarse problem solver in algebraic multigrid preconditioners.
It is part of the Trilinos software framework and designed to play a part in a robust and efficient solution of large-scale linear systems from real-world applications on current and emerging computers. 
In this paper, we discuss improvements to two sparse direct solvers, Basker and Tacho, and an algebraic preconditioner, FastILU, in the ShyLU-node package.
We present performance results with the sparse direct solvers for real application problems, namely, Basker for Xyce Circuit Simulations and Tacho for the Albany Land-Ice Simulation of Antarctica. FastILU has also been used in real-world applications, but in this paper, we illustrate its performance using 3D model problems.
\end{abstract}
\end{frontmatter}

%\keywords{
%performance portability, sparse linear algebra, dense linear algebra, graph algorithms.
%}

%\maketitle

\section{Introduction}

ShyLU-node is a software package designed to provide 
the scalable solution of the linear system of equations, $Ax=b$,
on a compute node (either on multicore CPUs or on a single GPU).
It is one of the open-source Trilinos software packages~\cite{trilinos} used to support large-scale scientific or engineering applications (in particular, for mission-critical applications at Sandia and other US National Laboratories). Trilinos is maintained and available on GitHub~\cite{Trilinos-website}.

Within the Trilinos software stack,
ShyLU-node is used as a stand-alone linear solver for a global problem,
as the local subdomain solver for domain decomposition (DD) preconditioner, or
as the coarse-problem solver in an algebraic multigrid preconditioner.
In many cases, the scalable on-node solvers can play an important role for achieving a high performance of the overall linear solvers and of the application simulations as they can be a significant part of the total simulation time.
For instance, even though DD methods provide an effective framework to construct a scalable preconditioner for solving a large-scale linear system on a distributed-memory computer, 
its performance may greatly depend on the performance of the local subdomain solvers (as we will show in Section~\ref{sec:tacho}).

In this paper, we discuss recent improvements made to two algebraic sparse direct solvers, which are provided as separate subpackages in ShyLU-node:
\begin{itemize}
\item Basker: a multi-threaded sparse direct solver, targeting linear problems for circuit simulations, and
\item Tacho: a multi-frontal sparse direct solver, based on Kokkos \cite{Kokkos} for performance portability over different node architectures (e.g., NVIDIA or AMD GPUs).
\end{itemize}
Besides these two sparse direct solvers, ShyLU-node also provides an algebraic preconditioning subpackage,
\begin{itemize}
\item FastILU: an iterative-variant of level-based incomplete LU (ILU) factorization and sparse-triangular solver, similar to the one proposed in~\cite{Chow:2015}. It is also based on Kokkos and provides the portable performance over different node architectures.
\end{itemize}
Although these subpackages have been introduced in their original papers~\cite{Tacho,Booth:2016}, in order to support real-world applications, their robustness, performance, and functionalities have been significantly enhanced over the years. To motivate the use of the solvers by many more users, in this paper, we summarize their current state and show the performance of the solvers for real-world applications, i.e., Basker for Xyce Circuit Simulations~\cite{doi:10.1142/9781860949630_0021,Xyce} and Tacho for Albany Land-Ice Simulation~\cite{gmd-11-3747-2018,Albany-webpage} of Antarctica.
While FastILU has been also used in real-world applications such as thermal simulation~\cite{aria}, to demonstrate its performance, in this paper, we use 3D model problems. To show the portability of the solvers, we conducted each experiments on different architectures
(ShyLU-Basker on Intel multi-core CPUs, Tacho on an NVIDIA A100 GPU, and FastILU on NVIDIA H100 GPU). % and AMD Mi250 GPU).

In the following sections, after listing related work (Section~\ref{sec:related}),
we first introduce Trilinos sofware framework (Section~\ref{sec:trilinos}) and describe ShyLU-node's three subpackages, Basker, Tacho, and FastILU, along with their performance results (Sections~\ref{sec:basker}, \ref{sec:tacho}, and \ref{sec:fastilu}, respectively). We then conclude with our final remarks (Section~\ref{sec:conclusion}).

\section{Related Work}\label{sec:related}

Though we only discuss the sparse direct solvers in ShyLU-node, there are several other on-node sparse direct solvers, and their descriptions can be found in survey papers~\cite{Duff:1984,Davis:2016}.
For instance, on the shared-memory CPUs, open-source SuperLU-MT package~\cite{SuperluMT} and Pardiso package~\cite{PARDISO} from Intel Math Kernel Library 
implement multi-threaded sparse direct solvers.
For a GPU, open-source Cholmod package~\cite{Cholmod} can offload the large dense blocks to a GPU for solving a sparse symmetric positive definite (SPD) problem, while the vendor-provided solvers such as CuSolver~\cite{cusolver} and RocSolver~\cite{rocsolver} implement a sparse direct solvers on their GPU. 

Ginkgo~\cite{ginkgo} is another software package that provides iterative variants of both level-set and threshold-based ILU on a GPU, while Kokkos-Kernels~\cite{KK} implements a standard level-set based ILU and an iterative variant of the threshold-based ILU using Kokkos.

This paper focuses on the recent improvements made to ShyLU-node, and the extensive comparison to other packages is out of its scope.

\section{Overview -- Trilinos}\label{sec:trilinos}

\begin{table*}[th]\small
\centerline{
  \subfloat[Linear Solver Packages\label{fig:trilinos-solvers}]{
  \begin{tabular}{ll}
  \multicolumn{2}{l}{Linear Solvers} \\
  \hline
  {\tt ShyLU}   & Distributed DD preconditioner (ShyLU-dd with {\tt FROSch}~\cite{frosch}) and\\
                & on-node factorization-based local solvers (ShyLU-node with {\tt Basker}~\cite{Booth:2016} and {\tt Tacho}~\cite{Tacho})\\
  {\tt Amesos2}~\cite{Amesos2} & Direct solver interfaces (e.g., {\tt KLU}~\cite{Davis:2010}, {\tt PaRDISO MKL}, {\tt SuperLU}~\cite{superlu}, {\tt MUMPS}~\cite{mumps}, {\tt Basker}, {\tt Tacho})\\
  {\tt Ifpack2}~\cite{Ifpack2} & Algebraic preconditioners (ILU, relaxation, one-level Schwarz)\\
  {\tt MuLue}~\cite{MueLu}     & Algebraic multigrid solver\\
  {\tt Teko}~\cite{Teko}       & Block preconditioner (for multi-physics problems)\\
  {\tt Belos}~\cite{Amesos2}   & Krylov solvers (e.g., CG, GMRES, BiCG, and
                                 their communication-avoiding or pipelined variants)\\
  \end{tabular}
  }
}
\centerline{
  \subfloat[Performance Portability Layers.\label{fig:trilinos-portability}]{
  \begin{tabular}{ll}  
  \multicolumn{2}{l}{Portable Performance} \\
  \hline
  {\tt Tpetra}~\cite{tpetra}     & Distributed sparse/dense matrix-vector operations \\
  {\tt Kokkos-Kernels}~\cite{KK} & Performance portable on-node graph and sparse/dense matrix operations\\
  {\tt Kokkos}~\cite{Kokkos}     & C++ programming model to provide performance portability on different node architectures\\
                                 &  (e.g., CPUs, NVIDIA/AMD GPUs)
  \end{tabular}
  }
}
\caption{Trilinos linear solver, and supporting, packages}\label{fig:trilinos}
\end{table*}

Trilinos~\cite{trilinos} is a collection of software packages used to support large-scale scientific and engineering applications. It is open-source software available on GitHub~\cite{Trilinos-website} with active developers and workflows to maintain the robustness of the overall software framework and to provide prompt support for its user communities~\cite{Milewicz:2022}. 

Table~\ref{fig:trilinos-solvers} shows the core Trilinos software packages for solving linear systems of equations, while Table~\ref{fig:trilinos-portability} lists the Trilinos software packages that enable our solvers to be portable to different computer architectures, using a single code base.
These solver packages can be combined to build a flexible and adaptable solver to address specific needs for solving a specific type of linear systems arising from an application. Trilinos is also C++ templated, for example, based on scalar types (double or float, and real or complex), local and global integer types (int or long long), and on-node programming models (serial, OpenMP, CUDA, HIP, etc.), providing further flexibility for its users. For instance, each matrix or vector, or preconditioner, can be of a different precision, enabling a mixed-precision operator or solver.

\begin{figure*}[t]
\scriptsize
\centerline{
\subfloat[Tacho Interface.\label{lst:tacho}]
{
  \input{codes/tacho}
}
\subfloat[Amesos2 Interface.\label{lst:amesos2}]
{
 \input{codes/amesos2}
}
}
\caption{Trilinos Linear Solver Interfaces.}
\end{figure*}

The Trilinos solvers consist of three steps. 
For instance, to solve the linear systems, $Ax=b$, based on the direct or approximate LU factorization of the input matrix $A$, ShyLU-node, in particular, performs the following tasks at each step: 
\begin{enumerate}
\item \emph{Symbolic Analysis} uses only 
      the sparsity structure of the coefficient matrix $A$.
      This step needs to be performed once for multiple solves with a fixed sparsity structure of $A$.
      The symbolic analysis often involves matrix reordering to expose more parallelism
      and to reduce the number of new nonzero entries, \emph{fills}, 
      introduced during the symbolic or numerical factorization of $A$. 
      It then analyzes and sets up the internal data structures for
      computing and storing the lower and upper triangular factors of the matrix $A$.
      
\item \emph{Numerical Setup}
      copies the numerical values of the coefficient matrix $A$
      into the internal data structures and computes
      the numerical factorization of the matrix, $P_r A P_c = LU$, 
      where $L$ and $U$ are the lower and upper triangular factors, and $P_r$ and $P_c$ are the
      row and column reordering computed during the symbolic analysis. To enhance the numerical stability of the factorization, additional scaling and permutation may be applied to the matrix $A$ during the numerical factorization as we will describe in the following sections.

\item \emph{Solve}
      computes the solution vectors $x$ based on forward and backward substitutions for a given set of right-hand-side vectors $b$. It performs sparse-triangular solve 
      using the lower and upper triangular factors, $L$ and $U$, computed by the numerical factorization phase.
\end{enumerate}
This software design is motivated by our application needs.
For instance, for many of our applications, the matrix sparsity structure stays the same throughout the entire simulation, and
only the numerical values of the matrix $A$ change, allowing us to amortize the cost of the \emph{symbolic analysis}, which is based only on the sparsity structure of the matrix.
Moreover, a sequence of the solutions often needs to be computed with different right-hand-side vectors but with the same matrix $A$,
allowing us to reuse the \emph{numerical factorization} of the matrix
for multiple \emph{solves}.
Hence, this software design provides the flexibility and improve the overall performance of the solver. 
The frequency and time spent for each phase depends on the application.
However, we typically focus on accelerating the performance of the numerical and solve phase, while we perform the symbolic factorization mostly on a single CPU core because the execution time for the symbolic factorization is often amortized over the whole simulation time (in addition, a sequential but high-quality analysis is often more important than a parallel analysis that may lower the quality, e.g., more fills).

Each software package is callable as an independent stand-alone solver. 
Nevertheless, to make it easier for the users to switch between different solvers,
Amesos2~\cite{Amesos2} and Ifpack2~\cite{Ifpack2} packages of Trilinos provide the uniform interface to different direct solvers and to algebraic preconditioners, respectively. For instance,  both Basker and Tacho are available through the Amesos2 interface, while FastILU is available through the Ifpack2 interface.
As an example, Listings~\ref{lst:tacho} and~\ref{lst:amesos2} show the Tacho and Amesos2 solver interfaces, respectively.

\section{Basker}\label{sec:basker}

Basker~\cite{Booth:2016} is a multi-threaded implementation of an algebraic sparse direct solver, KLU~\cite{Davis:2010}. It is designed, in particular, for the matrix sparsity structures arising from circuit simulations.
Unlike KLU or Basker, most existing sparse direct solvers are designed to take advantage of the \emph{supernodal} dense block structures, which are typical for the sparse matrices, arising from the mesh-based discretization of partial differential equations (PDE).
In contrast to these mesh-based matrices, the matrices from circuit simulations are sparser with heterogeneous sparsity structures and do not exhibit large enough dense blocks for the solvers to leverage for the performance gain (even after matrix reordering), or the supernodal approach could lead to a performance degradation.
In order to avoid the overhead associated with forming and operating on the supernodal blocks,
both KLU and Basker are based on a column-wise sparse LU factorization.

In addition, the circuit matrices are non-symmetric and some of them can be reordered into
a Block Triangular Form (BTF)~\cite{BTF} with small sparse diagonal blocks.
After the matrix is reordered into a BTF structure, KLU and Basker need to factorize only the diagonal blocks, reducing the computational costs of the factorization. Since these sparse diagonal blocks can be factorized independently, Basker uses multiple threads to factorize them in parallel (one thread per diagonal block).

\begin{figure*}
\centerline{
\subfloat[{\tt circuit\_4} after BTF.\label{fig:circut-4}]{
  \includegraphics[width=0.3\textwidth]{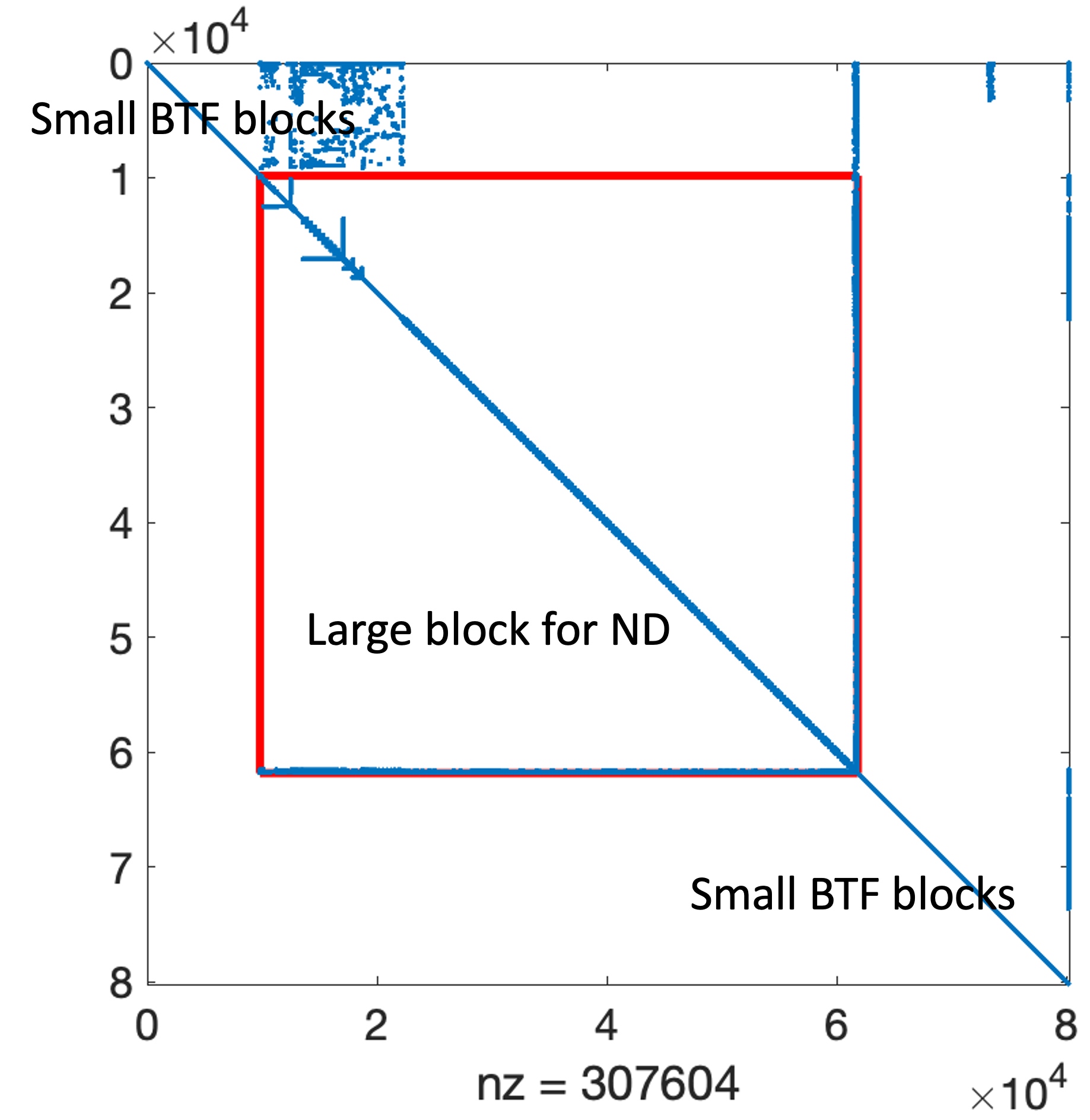}}
\quad\quad\quad\quad
\subfloat[Xyce Test Matrices (`n' is the size of $A$ and `nnz' is the number of non-zeroes in $A$).\label{tab:xyce}]{
\begin{tabular}{r|rrrr}\footnotesize
ID & n    & nnz/n   & \# BTF blocks & largest block size\\
\hline
1 & 682K & 5.7     & 707           & 681K\\
3 & 117K  & 4.1    & 79            & 117K\\
5 & 199K  & 12.6   & 11            & 199K\\
7 & 1.94M & 4.9    & 4105          & 1.94M\\
8 & 1.95M & 4.9    & 4105          & 1.94M\\
9 & 6.24M & 3.9    & 18            & 3.36M\\ %2.89M, \footnote{ShyLU-Basker currently supports only one big block.}\\
\end{tabular}
\vspace{1cm}}
}
\caption{BTF Structure and Properties of Test Matrices}
\end{figure*}

\begin{figure}[t]\footnotesize
\centerline{
\subfloat[Matrix ordering based on nested dissection.]{
  \includegraphics[width=0.25\textwidth]{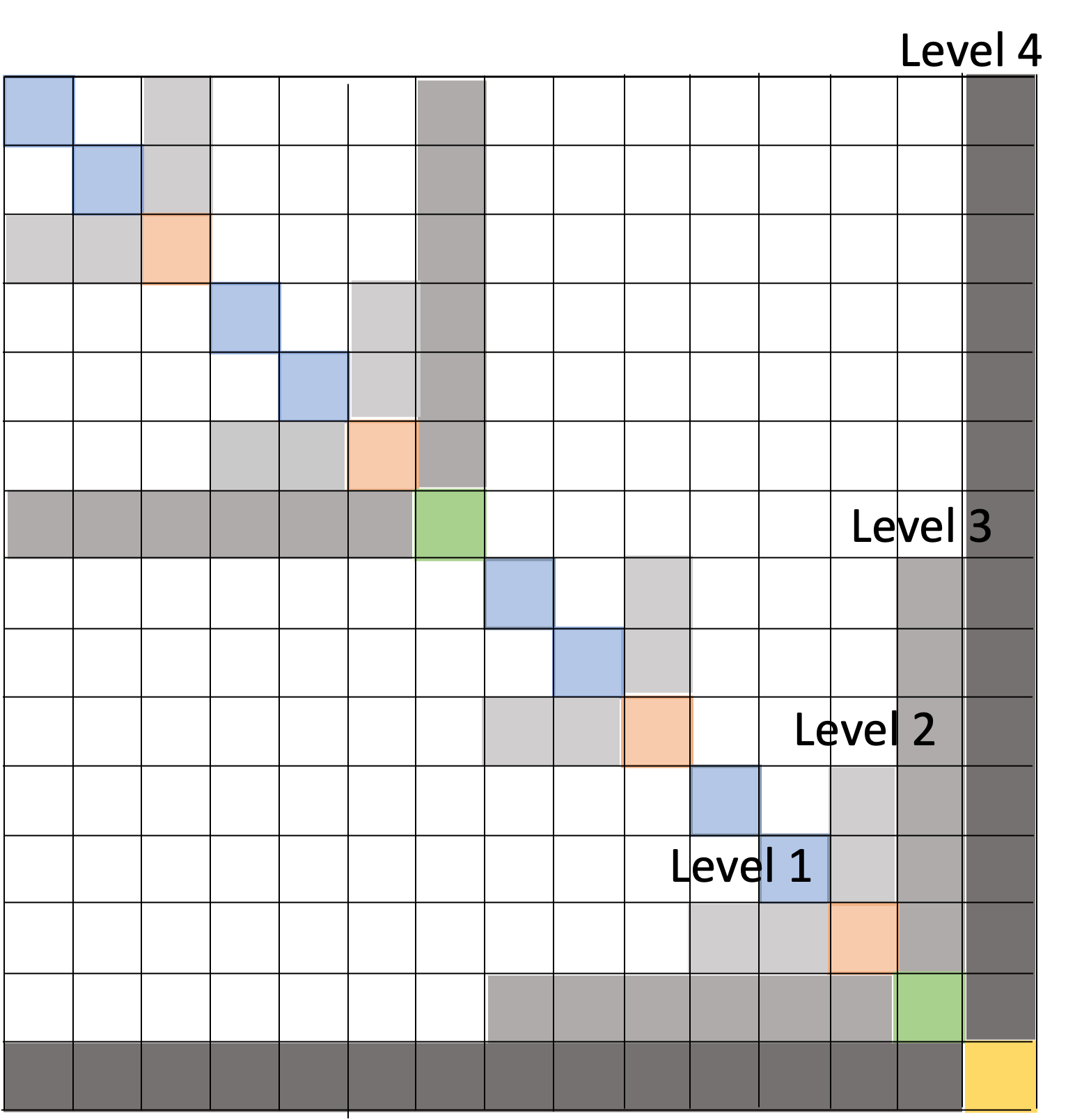}}
}
%\quad\quad\quad\quad
\centerline{
\subfloat[Parallel level-set scheduling of matrix factorization.]{
  \includegraphics[width=0.35\textwidth]{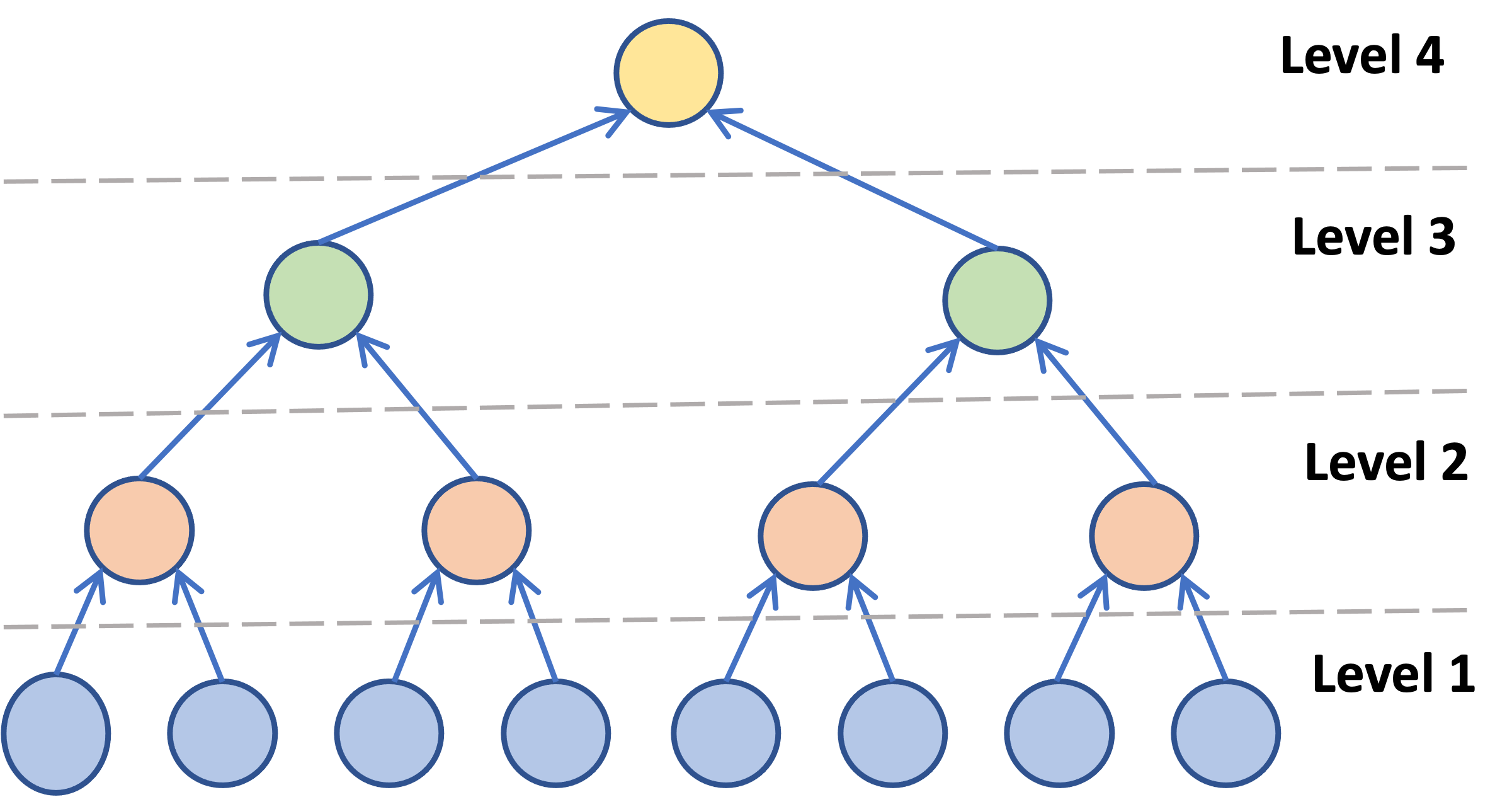}}
}
\caption{Nested dissection of sparse matrix (eight threads).}\label{fig:nd}
\end{figure}

\begin{figure*}[t]
\centerline{
\subfloat[Xyce Test 1 (BTF).]{
  \includegraphics[width=0.25\textwidth]{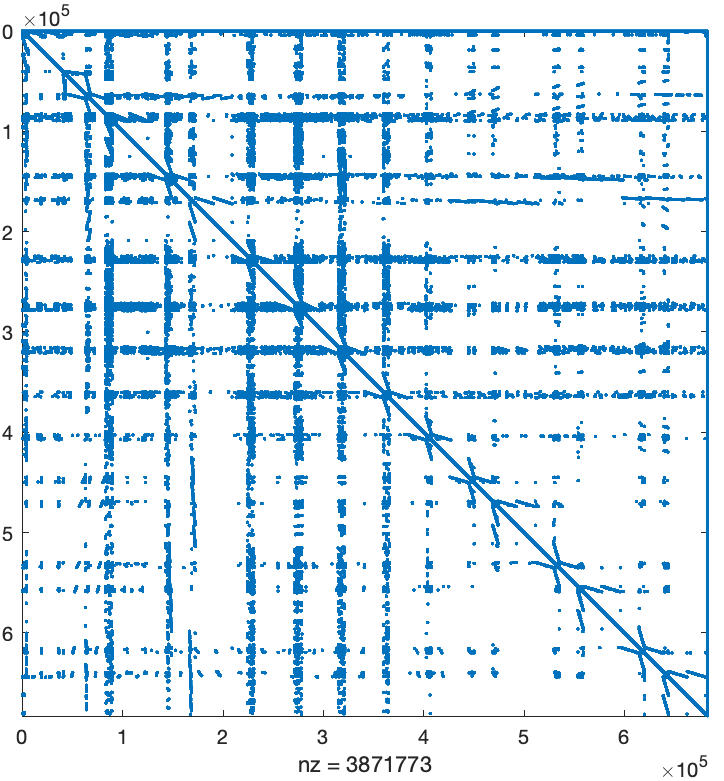}}
\quad\quad
\subfloat[Xyce Test 3 (BTF).]{
  \includegraphics[width=0.25\textwidth]{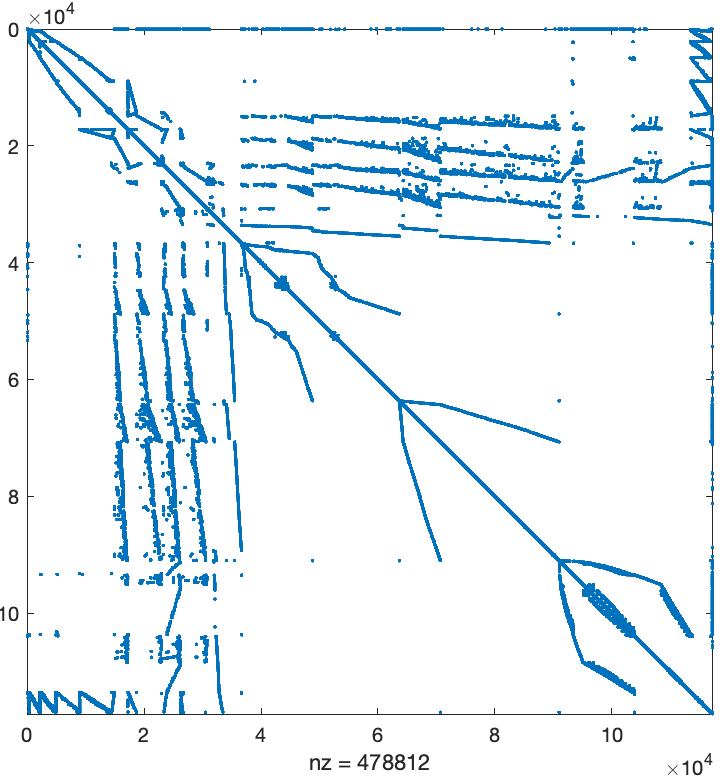}}
\quad\quad
\subfloat[Xyce Test 5 (BTF).]{
  \includegraphics[width=0.25\textwidth]{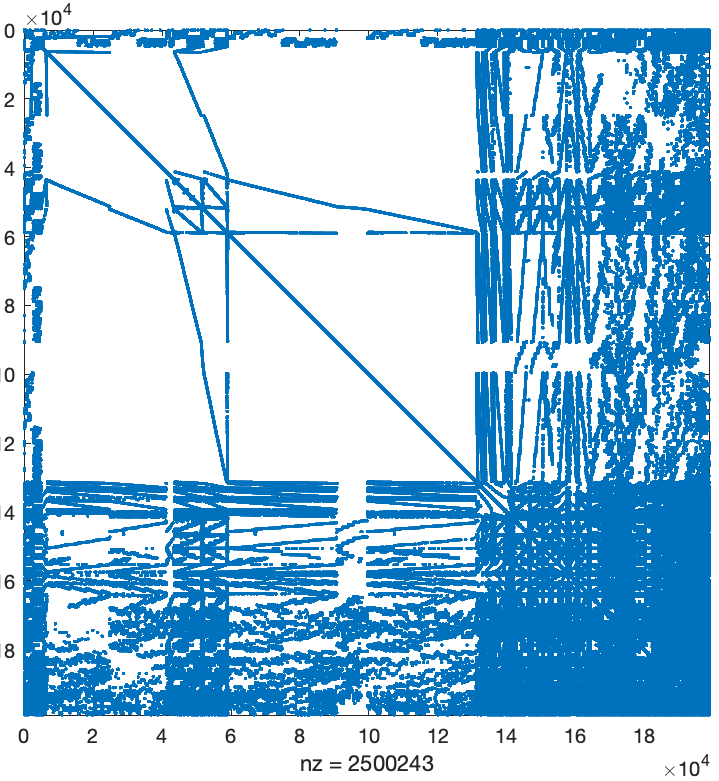}}
}
\centerline{
\subfloat[Xyce Test 1 (ND).]{
  \includegraphics[width=0.25\textwidth]{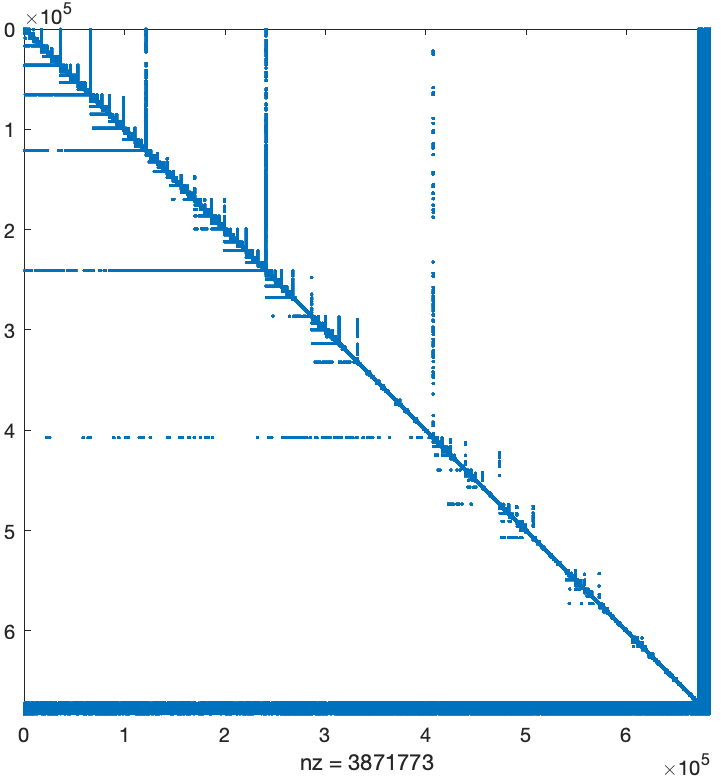}}
\quad\quad
\subfloat[Xyce Test 3 (ND).]{
  \includegraphics[width=0.25\textwidth]{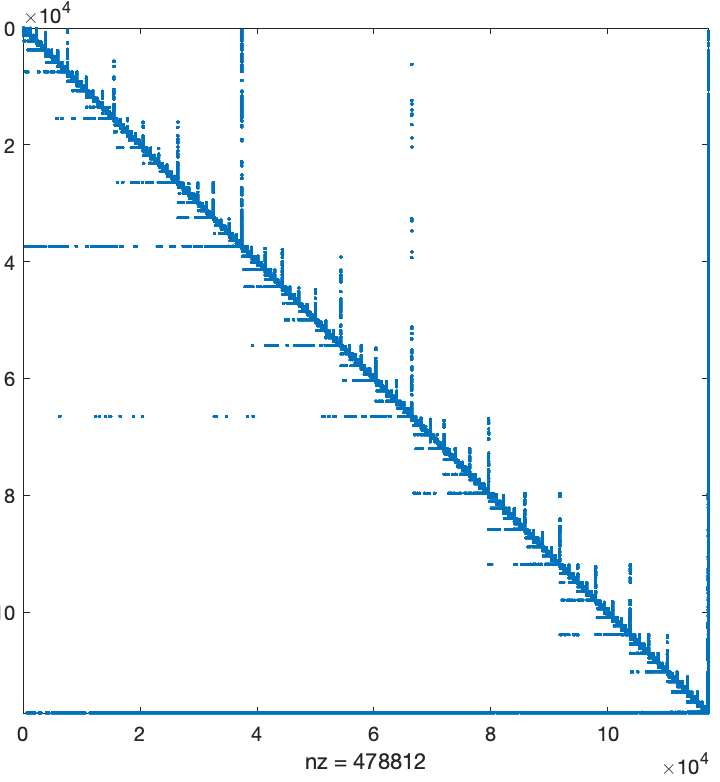}}
\quad\quad
\subfloat[Xyce Test 5 (ND).]{
  \includegraphics[width=0.25\textwidth]{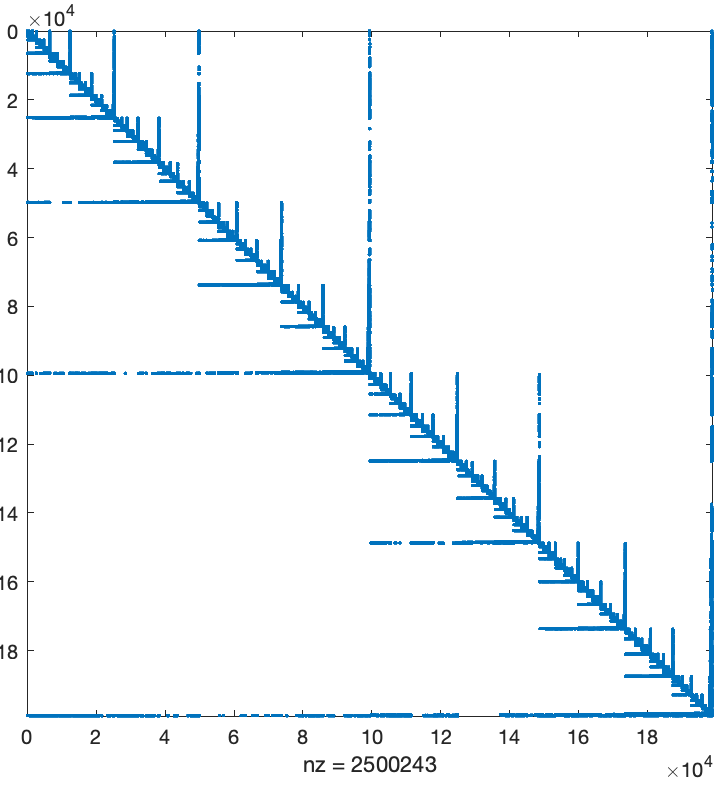}}
}
\caption{Sparsity Patterns of Xyce Test Matrices (see Table~\ref{tab:xyce} for the matrix descriptions).\label{fig:btf}}
\end{figure*}

As can be seen in Figure~\ref{fig:circut-4},
many of the matrices from modern circuit simulation do not exhibit the BTF structure (more symmetric in the structure) and contain large diagonal blocks after the BTF reordering. 
To improve thread scalability, after factorizing the small diagonal blocks in parallel,
Basker factorizes the largest diagonal block using multiple threads, based on the standard level-set scheduling scheme~\cite{Anderson:1989}.
To expose parallelism,
the large diagonal block is first reordered based on the nested dissection (ND) algorithm
of the METIS software package~\cite{metis} such that the number of interior leaf blocks is the same as the number of threads (see Figure~\ref{fig:nd} for an illustration).
Then, at each level, independent ND blocks (interior leaf blocks or interface separator blocks) are factorized in parallel. Since each block typically does not contain a large dense block, each ND block is factorized column by column.
For a large matrix, the leaf block can be relatively large. To reduce the number of fills in the LU factors, the leaf nodes are reordered using either the Approximate Minimum Degree (AMD) algorithm of Suite Sparse software~\cite{amd} or the ND algorithm of METIS,
where we found that METIS often leads to a longer symbolic setup time, but fewer fills and shorter numerical factorization time, especially for larger circuit matrices.
As Figure~\ref{fig:btf} shows, the large blocks can be often partitioned with small separator blocks, leading to good thread parallelization of factorizing the large block. Since the numerical factorization typically takes longer than each solve, Basker focuses on accelerating its numerical factorization time using threads while the sparse-triangular has not been parallelized.\footnote{The sparse-triangular solve of the large BTF block has recently been parallelized by solving the independent ND blocks in parallel at the leaf level.}

During the numerical factorization,
partial pivoting is applied only within the small diagonal BTF block or within each ND block of the large diagonal BTF block. Since swapping the pivot row with the diagonal row often increases the number of fills, pivoting is performed only when the magnitude of the diagonal entry is less than the magnitude of the largest entry in the column, multiplied by a threshold (i.e., $|a_{ii}| < \tau \max_k |a_{k,i}|$, where the default value of $\tau$ is $10^{-3}$). When the column does not contain any nonzero entries, Basker also provides an option to replace the diagonal entry with a small perturbation $\epsilon\|A_{ii}\|_1$, where $\|A_{ii}\|_1$ is the $1$-norm of the diagonal block and $\epsilon$ is the working machine epsilon. To recover the required solution accuracy, Amesos2 provides the option to perform iterative refinement with the computed LU factors.

\begin{table*}\footnotesize
\centerline{
\subfloat[ShyLU Basker using AMD reordering.]{
\begin{tabular}{r|ccccc}
   & \multicolumn{5}{c}{Number of Threads}\\
 ID         & 1                  & 2                   & 4                   & 8 & 16\\
 \hline\hline
 1 & 0.83 + 0.39 + 0.02 & 4.20 + 0.27 + 0.02 & 4.46 + 0.21 + 0.02 & 4.34 + 0.19 + 0.02 & 4.33 + 0.18 + 0.02\\
   & \multicolumn{1}{r}{\bf 0.41} & \multicolumn{1}{r}{\bf 0.29} & \multicolumn{1}{r}{\bf 0.23} & \multicolumn{1}{r}{\bf 0.21} & \multicolumn{1}{r}{\bf 0.20}\\
 \hline
 3 & 0.12 + 0.03 + 0.01 & 0.47 + 0.03 + 0.01 & 0.47 + 0.02 + 0.01& 0.48 + 0.02 + 0.01 & 0.49 + 0.02 + 0.01\\ 
   & \multicolumn{1}{r}{\bf 0.04} & \multicolumn{1}{r}{\bf 0.04} & \multicolumn{1}{r}{\bf 0.03} & \multicolumn{1}{r}{\bf 0.03} & \multicolumn{1}{r}{\bf 0.03}\\
 \hline
 5 & 0.83 + 1.20 + 0.01 & 2.55 + 0.69 + 0.01 & 2.57 + 0.40 + 0.01 & 2.57 + 0.23 + 0.01 & 2.58 + 0.17 + 0.01\\
   & \multicolumn{1}{r}{1.21} & \multicolumn{1}{r}{0.70} & \multicolumn{1}{r}{0.41} & \multicolumn{1}{r}{0.24} & \multicolumn{1}{r}{0.18}\\
 \hline
 7 & 3.71 + 8.94 + 0.10 & 16.1 + 4.69 + 0.12 & 16.0 + 3.33 + 0.10 & 16.8 + 2.13 + 0.10 & 16.3 + 1.05 + 0.10\\
   & \multicolumn{1}{r}{9.04} & \multicolumn{1}{r}{4.81} & \multicolumn{1}{r}{3.43} & \multicolumn{1}{r}{2.23} & \multicolumn{1}{r}{1.15}\\
 \hline
 8 & 13.4 + 9.37 + 0.10 & 16.1 + 5.01 + 0.10 & 16.2 + 3.38 + 0.10 &  16.9 + 2.15 + 0.10 & 17.9 + 1.13 + 0.10\\
   & \multicolumn{1}{r}{9.47} & \multicolumn{1}{r}{5.11} & \multicolumn{1}{r}{3.48} & \multicolumn{1}{r}{2.25} & \multicolumn{1}{r}{1.23}\\
 \hline
 9 & 5.37 + 1.10 + 0.22 & 21.8 + 1.18 + 0.22 & 21.8 + 1.04 + 0.20 & 21.8 + 1.02 + 0.21 & 22.3 + 1.10 + 0.32\\
   & \multicolumn{1}{r}{\bf 1.32} & \multicolumn{1}{r}{\bf 1.40} & \multicolumn{1}{r}{\bf 1.24} & \multicolumn{1}{r}{\bf 1.23} & \multicolumn{1}{r}{\bf 1.42}\\
% \hline
% 14 & 19.2 + 309.1 + 1.10 & 
\end{tabular}}}
\centerline{
\subfloat[ShyLU Basker using METIS reordering.]{
\begin{tabular}{r|ccccc}
            & 1                  & 2                   & 4                   & 8 & 16\\
 \hline\hline
 1 & 7.33 + 0.38 + 0.03 & 8.53 + 0.26 + 0.03 & 7.38 + 0.21 + 0.03 & 6.68 + 0.20 + 0.03 & 6.00 + 0.19 + 0.03\\
   & \multicolumn{1}{r}{\bf 0.41} & \multicolumn{1}{r}{\bf 0.29} & \multicolumn{1}{r}{0.24} & \multicolumn{1}{r}{0.23} & \multicolumn{1}{r}{0.22}\\
 \hline
 3 & 0.52 + 0.04 + 0.01 & 0.80 + 0.03 + 0.01 & 0.78 + 0.02 + 0.01 & 0.75 + 0.02 + 0.01 & 0.71 + 0.02 + 0.01\\ 
   & \multicolumn{1}{r}{0.05} & \multicolumn{1}{r}{\bf 0.04} & \multicolumn{1}{r}{\bf 0.03} & \multicolumn{1}{r}{\bf 0.03} & \multicolumn{1}{r}{\bf 0.03}\\
 \hline
 5 & 6.43 + 0.87 + 0.02 & 6.27 + 0.48 + 0.02 & 5.30 + 0.29 + 0.02 & 4.72 + 0.19 + 0.02 & 4.02 + 0.15 + 0.02\\
   & \multicolumn{1}{r}{0.89} & \multicolumn{1}{r}{0.50} & \multicolumn{1}{r}{0.31} & \multicolumn{1}{r}{0.21} & \multicolumn{1}{r}{\bf 0.17}\\
 \hline
 7 & 20.7 + 2.43 + 0.11 & 40.1 + 1.46 + 0.12 & 39.3 + 0.92 + 0.13 & 35.6 + 0.80 + 0.13 & 33.7 + 0.64 + 0.13\\
   & \multicolumn{1}{r}{2.54} & \multicolumn{1}{r}{\bf 1.58} & \multicolumn{1}{r}{\bf 1.05} & \multicolumn{1}{r}{\bf 0.93} & \multicolumn{1}{r}{\bf 0.77}\\
 \hline
 8 & 21.4 + 2.50 + 0.12 & 40.0 + 1.43 + 0.12 & 41.8 + 1.11 + 0.12 & 37.4 + 0.82 + 0.13 & 34.0 + 0.61 + 0.13\\
   & \multicolumn{1}{r}{2.62} & \multicolumn{1}{r}{\bf 1.55} & \multicolumn{1}{r}{\bf 1.23} & \multicolumn{1}{r}{\bf 0.95} & \multicolumn{1}{r}{\bf 0.74}\\
 \hline
 9 & 5.98 + 1.12 + 0.39 & 35.8 + 1.35 + 0.40 & 34.7 + 1.11 + 0.40 & 33.9 + 1.07 + 0.42 & 33.4 + 1.14 + 0.42\\
   & \multicolumn{1}{r}{1.51} & \multicolumn{1}{r}{1.75} & \multicolumn{1}{r}{1.51} & \multicolumn{1}{r}{1.49} & \multicolumn{1}{r}{1.56}
% \hline
% 14 & 56.4 + 30.4 + 0.99 & 123.4 + 21.4 + 1.03 & 138.5 + 14.4 + 1.03 &  
\end{tabular}}}
\centerline{
\subfloat[Pardiso MKL.]{
\begin{tabular}{r|ccccc}
            & 1 & 2 & 4 & 8 & 16\\
 \hline\hline
 1 & 2.78 + 9.94 + 0.06 & 2.71 + 6.06 + 0.04 & 2.63 + 3.90 + 0.03 & 2.67 + 2.36 + 0.03 & 2.69 + 1.54 + 0.03\\
   & \multicolumn{1}{r}{10.00} & \multicolumn{1}{r}{6.10} & \multicolumn{1}{r}{3.93} & \multicolumn{1}{r}{2.39} & \multicolumn{1}{r}{1.57}\\
 \hline
 3 & 0.36 + 0.10 + 0.01 & 0.37 + 0.08 + 0.01 & 0.36 + 0.08 + 0.01 & 0.36 + 0.06 + 0.01 & 0.36 + 0.06 + 0.01\\
   & \multicolumn{1}{r}{0.11} & \multicolumn{1}{r}{0.09} & \multicolumn{1}{r}{0.09} & \multicolumn{1}{r}{0.07} & \multicolumn{1}{r}{0.07}\\
 \hline
 5 & 1.42 + 0.43 + 0.02 & 1.39 + 0.30 + 0.01 & 1.36 + 0.19 + 0.01 & 1.34 + 0.13 + 0.02 & 1.36 + 0.20 + 0.01\\
   & \multicolumn{1}{r}{\bf 0.45} & \multicolumn{1}{r}{\bf 0.31} & \multicolumn{1}{r}{\bf 0.20} & \multicolumn{1}{r}{\bf 0.15} & \multicolumn{1}{r}{0.21}\\
 \hline
 7 & 7.10 + 2.18 + 0.32 & 6.91 + 1.62 + 0.22 & 6.79 + 1.12 + 0.16 & 6.73 + 0.96 + 0.14 & 6.84 + 0.75 + 0.13\\
   & \multicolumn{1}{r}{\bf 2.50} & \multicolumn{1}{r}{1.84} & \multicolumn{1}{r}{1.28} & \multicolumn{1}{r}{1.10} & \multicolumn{1}{r}{0.88}\\
 \hline
 8 & 7.88 + 2.20 + 0.33 & 7.46 + 1.92 + 0.22 & 7.52 + 1.34 + 0.16 & 7.49 + 1.13 + 0.14 & 7.31 + 1.83 + 0.13\\
   & \multicolumn{1}{r}{\bf 2.53} & \multicolumn{1}{r}{2.14} & \multicolumn{1}{r}{1.50} & \multicolumn{1}{r}{1.27} & \multicolumn{1}{r}{1.96}\\
 \hline
 9 & 20.8 + 5.16 + 0.28 & 20.3 + 4.21 + 0.22 & 19.9 + 3.45 + 0.19 & 19.8 + 2.73 + 0.18 & 20.3 + 2.20 + 0.20\\
   & \multicolumn{1}{r}{5.44} & \multicolumn{1}{r}{4.43} & \multicolumn{1}{r}{3.64} & \multicolumn{1}{r}{2.91} & \multicolumn{1}{r}{2.40}
 %\hline
 %14 & 45.9 + 13.4 + 1.87 & 
\end{tabular}}}
\caption{Performance of ShyLU-Basker and Pardiso MKL (the first row shows symbolic + numeric + solve time in seconds, while the second row shows the total of the numeric and solve time)
for the test matrices show in Table~\ref{tab:xyce}.\label{tab:xyce-perf}}
\end{table*}

Basker provides several solver options including:
\begin{itemize}
\item {\bf Maximum Cardinality Matching} of input matrix or each diagonal block: This is used to reorder the matrix such that the diagonal entry has a self-edge. Although the matrix may explicitly store zero entries (which may become non-zero during the simulation), this will increase the likelihood of the diagonal entries being non-zero, enhancing the stability and performance of the numerical factorization. Since this option relies only on the sparsity structure of the matrix, it is computed during the symbolic analysis, and can be reused for multiple numerical factorizations.

\item {\bf Maximum Weighted Matching} of input matrix or each diagonal block: This is used to move the large nonzero entries to the diagonal. Since pivoting is applied only within each diagonal block, the stability of the overall solver may be improved by applying this option to the input matrix. In addition, this option could improve the performance because many of the diagonal entries could be used as pivots, reducing the number of row interchanges and of fills. Since the numerical values of the matrix are used to compute the matrix ordering, this option is part of the numerical factorization.
%\ichi{``refactorization'' by reusing the matching from a seed matrix?}

\item {\bf Dynamic Reallocation}: Since partial pivoting is used within each diagonal block, the number of fills in the LU factors may increase. As a result, the storage allocated by the symbolic analysis may not be sufficient to store the LU factors. When a larger storage is needed, Basker reallocates the storage and restarts the factorization.
\end{itemize}
Since the first publication~\cite{Booth:2016}, both robustness and performance of Basker have been improved and it is one of the primary solvers used for Xyce circuit simulations~\cite{Keiter:2011,Xyce}.

To ensure the robust solutions of the linear problems for the whole simulation steps of different Xyce circuit simulations~\cite{Xyce}, the direct methods are often used. In addition, due to the heterogeneous structures of the matrices, the distributed-memory sparse direct solver can suffer from a high symbolic analysis cost or from the large memory requirement. Hence, though the Xyce simulation is typically ran using multiple MPI processes (to accommodate the memory needs), to solve the linear system, in many cases, the coefficient matrix $A$ is gathered to a single MPI process and a shared-memory sparse direct solver, like Basker, is used. The rest of the Xyce simulation typically scales well using multiple MPI processes; hence, the linear solver often becomes the performance bottleneck of the simulation.

Table~\ref{tab:xyce} compares the performance of ShyLU-Basker with a vendor-optimized PARallel DIrect SOlver (PARDISO) from Intel Math Kernel Library (MKL). 
The experiments were carried out on Intel Xeon Platinum CPUs, where the code was compiled using Intel 2021.5.0 compilers.
Throughout the circuit simulation, the matrix structure remains the same, but the numerical values of $A$ change at each step. Hence, the numerical factorization is recomputed to solve the linear problem at each simulation step, while the symbolic factorization needs to be performed only once and its cost is amortized over the simulation (e.g, thousands of steps). 
Therefore, we focus on the time required for the numerical factorization and solve phases. These test matrices are from modern circuit simulations and do not exhibit any BTF structure (see Figure~\ref{fig:btf}).
We observe that even without the BTF structure, ShyLU-Basker may outperform the vendor-optimized sparse direct solver, PARDISO MKL, depending of the sparsity structure of the matrix.

\section{Tacho}\label{sec:tacho}

Tacho~\cite{Tacho} implements a sparse direct solver based on multi-frontal method. 
Originally, Tacho only supported Cholesky factorization of a symmetric positive definite (SPD)
matrix, based on dynamic task scheduling provided through Kokkos on an NVIDIA GPU. Since then, both the capability and the portability of Tacho have been extended.
For instance, the current Tacho
is still based on Kokkos, but in order to provide performance portability to different node architectures, it now uses the
standard level-set scheduling for both numerical factorization and sparse-triangular solve (SpTRSV). Moreover, besides the Cholesky factorization of an SPD matrix, Tacho now supports various symmetric factorizations, namely, LDL$^T$ factorization of a symmetric indefinite matrix and LU factorization of a general matrix with a symmetric sparsity structure\footnote{Tacho now also supports sequential LDL$^T$ factorization of a skewed symmetric matrix.}.

Tacho also includes new solver options.
For instance, in many of our applications, multiple linear systems with the same matrix
but with different right-hand-side vectors need to be solved in sequence. To accommodate such needs,
Tacho provides four different variants of SpTRSV with trade-offs between the numerical stability, numerical setup cost, and triangular-solve performance\footnote{Similar variants of sparse-triangular solvers are implemented in Kokkos-Kernels~\cite{KK}, and their performance has been compared in~\cite{Yamazaki:2020}.}.
Going from variant 1 to variant 4, these four solver variants typically achieve higher performance of the sparse-triangular solve, but with the higher setup costs and potential for numerical instability. Our discussion below assumes a single right-hand-side vector, but Tacho supports solving multiple right-hand side vectors at once (e.g., using TRSM, GEMM, and SpMM instead of using TRSV, GEMV, and SpMV).
\begin{enumerate}
    \item [1] The first variant implements SpTRSV based on the level-set scheduling~\cite{Anderson:1989} of the
              dense supernodal blocks. At the lower level of the level-set scheduling tree, we typically have many small dense blocks, while we have a fewer but larger blocks at the higher level of the tree. Hence, at a lower level, Tacho launches batched BLAS kernels to perform the operations on the small blocks in parallel (TRSV to compute the solution block with the diagonal block, $x_i := L_{i,i}^{-1} b_i$, and then GEMV to update the remaining solution block $b_k := b_k - L_{k,i}x_i$), but then switches to calling a threaded BLAS kernels (from a vendor-optimized library such as CuBLAS or rocBLAS) at a higher level, where GPU streams are used to launch the threaded BLAS on different blocks in parallel.
    \item [2] Instead of calling a dense triangular solve (TRSV) with the diagonal blocks, this variant
              explicitly computes the inverse of the diagonal blocks and calls a dense matrix-vectors mulltiply (GEMV), which often leads to higher performance than TRSV. Since the inverses
              of the diagonal blocks are explicitly computed, the accuracy of the solution could deteriorate, depending on the condition number of the diagonal blocks.
    \item [3] In addition to inverting the diagonal blocks, this variant applies the diagonal inverse to 
              the corresponding off-diagonal blocks. This allows Tacho to merge the two GEMVs (one with the diagonal blocks and the other with the off-diagonal blocks) into a single GEMV at each level, reducing the number of kernel launches (where SpTRSV performance is often limited by either latency or bandwidth). 
              
              This algorithm is equivalent to the partitioned inverse~\cite{Alvarado:1993} based on the level-set partitioning of the supernodal block columns~\cite{Yamazaki:2020}, and algorithmically converts the SpTRSV to a sequence of Sparse-Matrix Vectors multiplies (SpMVs),
              \[
                 L^{-1} = \prod_{\ell=1}^{n_{\ell}} L^{-1}_{\ell},
              \]
              where $n_{\ell}$ is the number of levels, and $L_{\ell}$ is an identity matrix except that the supernodal columns, which belong to the $\ell$th level, are replaced with the corresponding columns in $L_{\ell}$.
              
              This method is often referred to as a no-fill method because the new fills introduced by the matrix inversion are restricted only within the non-empty blocks, which were already present in the computed LU factors (leading to no fills if the sizes of all the blocks are one). 
              
              Since the matrix inverse is computed at each level, the stability of the algorithm depends on the condition number of the partitioned matrix $L_{\ell}$ at each level.
    \item [4] The last variant is numerically the same as the third variant, but the inverse of the partitioned matrix at each level, $L_{\ell}^{-1}$,
              is stored in the sparse-matrix storage (Compressed Sparse Row) format. This allows Tacho to call a single Sparse-Matrix Multiple Vectors multiply (SpMV) kerenel at each level (without synchronizations between the levels), and also potentially reduces the number of nonzero entries that need to be stored in each block. On a GPU, Tacho uses the SpMV implementation from the vendor-optimized libraries (cuSparse or rocSparse).
\end{enumerate}
Tacho has been used as the local subdomain solver for various large-scale applications on GPU clusters,
including Structural Dynamics~\cite{sierra-sd},
Shape Optimization~\cite{plato}, and
Compressible Computational Fluid Dynamics~\cite{cfd}, and these solver options provide flexibility to achieve the overall high performance of the simulation with their specific needs.
For example, for the structural dynamics simulation, where the numerical factorization needs to be computed only once, the Variant 4 of SpTRSV is used, while as we show below, the Variant 1 is used for the Albany Land Ice simulation, where the numerical factors are recomputed for each solve.

\begin{table*}\small
%\centerline{
%\subfloat[CPU runs with SuperLU. \ichi{remove}]{
%\begin{tabular}{l|rrrrr}
%\# MPIs  & 4           & 8           & 16          & 32          & 64\\
%\hline
%Numeric & 6,111.4 (5,984.0) & 2,564.4 (2,429.9) & 1,059.1 (1,035.4) & 475.5 (442.7) & 151.6 (141.0)\\
%Solve   & 340.5 (115)       & 216.4 (120)       & 147.5 (152)       & 88.0 (154) & 38.9 (155) 
%\end{tabular}}}
%
\centerline{
\subfloat[CPU runs with Tacho.]{
\begin{tabular}{l|rrrrr}
\# MPI ranks  & 4           & 8           & 16          & 32          & 64\\
\hline
Numeric & 1,501.0 (1,455.6) & 624.6 (607.2)     & 300.4 (291.9)     & 170.6 (164.9) & 49.1 (44.9)\\
Solve   & 123.8 (115)       & 112.8 (120)       & 107.3 (152)       & 52.8 (154)    & 25.1 (155) 
\end{tabular}}}
\centerline{
\subfloat[GPU runs with Tacho.]{
\begin{tabular}{l|rrrrr}
\# MPI ranks  & 4           & 8           & 16          & 32          & 64\\
\hline
%Numeric  & 77.5 (69.1) & 48.7 (39.8) & 31.8 (23.2) & 26.7 (16.2) & 24.2 (10.0)\\
%Solve    & 14.4 (155)  & 9.8 (120)   & 8.6 (152)   & 7.4 (154)   & 6.1 (155)
%Numeric  & 86.2 (75.5) & 56.6 (46.3) & 41.2 (30.8) & 43.9 (25.3) & -- (--)\\
%Solve    & 14.1 (155)  & 9.6 (120)   & 8.6 (152)   & 7.6 (154)   & -- (--)
Numeric  & 77.5 (67.7) & 48.7 (29.0) & 31.8 (23.2) & 26.2 (16.6) & 21.8 (10.5)\\
Solve    & 14.4 (155)  & 9.4 (120)   & 8.6 (152)   & 7.7 (154)   & 7.1 (155)
\end{tabular}}}
\caption{Strong-scaling results on one Perlmutter GPU node ($2.2\times10^6$ DoFs).
The table shows the total time spent in seconds for the numerical setup and solve for the two-level DD solver, while the numbers in the parenthesis show the time spent in Tacho for the numerical factorization and the total number of iterations for the solve. 
\label{tab:albany-strong}}
\end{table*}

%Other solver options:
%\begin{itemize}
%\item 
%Graph Condensation: 
Another option provided by Tacho is to perform its symbolic analysis based on a ``condensed'' graph, which is generated by compressing the adjacency graph of the original matrix, for example, based on the underlying block structure for a multi-physics simulation. The original motivation for this option was to reduce the size of the graph and the cost of symbolic analysis (which is mostly sequential on a host CPU). However, it can also be used as a heuristic to enhance the numerical stability of the factorization. For instance, it has been used to ensure certain degree-of-freedoms (DoFs) to stay together in the supernodal block generated during the symbolic analysis, e.g., the multiple physics DoFs defined on each grid point or pivoting rows.  Since Tacho only pivots within each diagonal block, the numerical stability may be improved when these DoFs remain within the same diagonal block (this option is currently used for some of the multi-physics Albany Icesheet benchmark simulations to enhance the stability and for the skew-symmetric factorization to form the 2-by-2 pivots). A user could provide a mapping of the DoFs in the original matrix to the DoFs in the condensed matrix, or Tacho also provides an option to internally construct the condensed graph with a fixed block size when the DoFs, which belong to the same block, are stored contiguously in the sparse matrix storage (Compressed Sparse Row) format.

%\item Diagonal Perturbation
%\end{itemize}

As Tacho is still based on Kokkos, it is portable to different GPU architectures (NVIDIA or AMD) and also supports multi-threading on a shared-memory CPUs through OpenMP.
To demonstrate the performance of Tacho,
we show our experimental results where Tacho was used as the local subdomain solver for a multi-level overlapping additive Schwarz preconditioner (implemented in the ShyLU-DD software package called FROSch~\cite{frosch}).
This DD preconditioner (using overlap of one) was then combined with a Krylov subspace iterative solver, GMRES~\cite{Saad:1986}, to solve the nonsymmetric linear systems
arising during the Antarctica Icesheet simulation by Albany~\cite{albany}.

For this experiment,
we used the same meshes that were also used for the performance studies of Albany in~\cite{mali}. This performance benchmark runs for eight nonlinear solve steps and requires one numerical factorization for solving one linear system at each nonlinear step.
We considered GMRES to have converged when the relative residual
norm was reduced by six orders of magnitude and GMRES iteration was restarted after every 200 iterations (GMRES converged before the restart).
To improve the performance of local solve, Tacho reordered the local subdomain matrix using METIS before the factorization.
Since we need to recompute the numerical factorization at each step of the simulation, we used the variant 1 of SpTRSV with Tacho.

We compare the numerical setup and GMRES iteration time for Albany built for CPUs and GPUs on the Perlmutter Supercomputer at National Energy Research Scientific Computing (NERSC) Center. 
Each GPU node of Perlmutter has 64-core AMD EPYC 7763 (Milan) CPU and 4 NVIDIA A100 GPUs. 
The code was compiled using Cray's compiler wrapper 
with Cray LibSci version 23.2, CUDA version 11.7
and Cray MPICH version 8.1.
The GPU-aware MPI was not available on Perlmutter, and hence, all the MPI communications were performed through the CPU.\footnote{On other machines such as the Frontier supercomputer at the Oak Ridge Leadership Computing Facility, a significant GPU performance gain has been reported using GPU-aware MPI.}
We configured Trilinos such that all the local dense and sparse matrix operations are performed using CUBLAS and CuSparse, respectively.

The complexity of the local sparse direct solver scales more than linearly (e.g., $\mathcal{O}(n_i^6)$ computational and $\mathcal{O}(n_i^2)$ storage complexities for the local structured 3D problem of dimension $n_i^3$, using a nested dissection ordering).  As a result, for our CPU runs, FROSch obtained the best performance running one MPI per CPU core~\cite{Yamazaki:2023} (assuming that the coarse-space problem keeps the iteration counts constant with an increasing number of subdomains and that the coarse-problem solve or the MPI communication does not become significant in the total solver time).
For our GPU runs to have the same number of subdomains as our CPU runs,
we used the NVIDIA's Multi-Process Service (MPS) to launch multiple MPI processes on each A100 GPU (up to 16 MPI processes on each GPU).

\ignore{
\begin{table}\small
\centerline{
\subfloat[CPU runs with Tacho (64 MPI ranks / node).]{
\begin{tabular}{l|rrrrr}
\# nodes  & 1               & 4    & 16 \\
%\# DoFs   & $2.2\times10^6$ & $8.8\times10^6$\\
\hline
Numeric & 49.1 (44.9) & 53.0 (47.8) & 65.8 (59.3)\\
Solve   & 25.1 (155)  & 28.4 (160)  & 48.5 (253)
\end{tabular}}
}
%
%\quad\quad\quad
%
\centerline{
\subfloat[GPU runs with Tacho (32 MPI ranks / node).]{
\begin{tabular}{l|rrrrr}
\# nodes  & 1 & 4 & 16 \\
\hline
Numeric  & 26.2 (16.6) & 32.3 (19.6) \\
Solve    & 7.7 (154)   & 8.9 (167)
\end{tabular}}
}
\caption{Weak-scaling results on multiple Perlmutter GPU nodes (with about $2.2\times10^6$ DoFs / node),
where all the CPU cores were used for the CPU runs (1 MPI / core with 64 cores / node), while 8 MPI processes were launched on each GPU using MPS for our GPU runs (8 MPI ranks / GPU with 4 GPUs / node).\label{tab:albany-weak}}
\end{table}
}

\begin{table}\small
\centerline{
\begin{tabular}{l|rr|rr}
          & \multicolumn{2}{c|}{CPU runs } & \multicolumn{2}{c}{GPU runs }\\
\# nodes  & 1               & 4     & 1 & 4\\
\hline
Numeric & 49.1 (44.9) & 53.0 (47.8) & 26.2 (16.6) & 32.3 (19.6)\\
Solve   & 25.1 (155)  & 28.4 (160)  & 7.7 (154)   & 8.9 (167)
\end{tabular}
}
\caption{Weak-scaling results on multiple Perlmutter GPU nodes (with about $2.2\times10^6$ DoFs / node),
where all the CPU cores were used for the CPU runs (1 MPI / core with 64 cores / node), while 8 MPI processes were launched on each GPU using MPS for our GPU runs (8 MPI ranks / GPU with 4 GPUs / node). Tacho is used as the local subdomain solver on both CPUs and GPU.
The table shows the total time spent in seconds for the numerical setup and solve for the two-level DD solver, 
while the numbers in the parenthesis show the time spent in Tacho for the numerical factorization and the total number of iterations for the solve. \label{tab:albany-weak}}
\end{table}

Table~\ref{tab:albany-strong} shows the parallel strong-scaling performance results on one node of Perlmutter. The performance of the overall domain-decomposition solver depends critically on MPS. Nevertheless, when MPS can efficiently run multiple MPI processes on each GPU, we see that Tacho was able to obtain GPU acceleration for the local subdomain solver and for the overall domain-decomposition solver with the speedups of 17.7$\times$ and 2.6$\times$ using 4 and 64 MPI processes, respectively. In particular, for the last column of the table, our CPU experiment used all the CPU cores available on the GPU node; hence, our comparison is based on the CUP experiments using 16 CPU cores for each GPU used for the GPU experiment.

We have also used SuperLU 5.2.1 as the local subdomain solver on the CPU (compiled with {\tt cc} Cray compiler wrapper and linked to Cray LibSci library for its BLAS and LAPACK, while using METIS for reordering the matrix). For this particular simulation, Tacho, running with one thread, was more efficient than SuperLU. This could be because, while SuperLU performs the partial pivoting to maintain the numerical stability of the factorization, Tacho looks for the pivots only within each diagonal blocks for performance. For this particular simulation, Tacho was stable enough.

Similarly, Table~\ref{tab:albany-weak} shows the parallel weak-scaling results. Though the overall performance of the domain-decomposition solver may improve by having more MPI processes on a single GPU with MPS, there could be a significant overhead (e.g., GPU memory). Hence, for our weak-scaling experiments, we launched only 32 MPI processes on each node for our GPU runs (i.e., 8 MPI ranks / GPU), while 64 MPI processes were used on each node for our CPU runs (i.e., one MPI / core). %Hence on 1 and 4 compute nodes, we have 512 and 1024 MPI processes for our GPU and CPU run, respectively. We observed that the GPU acceleration was maintained over multiple compute nodes.

\section{FastILU}\label{sec:fastilu}

Preconditioners based on Incomplete LU (ILU) factorization are effective for some applications. Several parallel ILU factorization algorithms (e.g., based on level-set scheduling~\cite{Anderson:1989}) have been developed, along with a matrix reordering technique to enhance the parallelism~\cite{Jones:1994}. However, the performance of computing the ILU preconditioner, and the sparse-triangular solve (SpTRSV) needed to apply the preconditioner, depends on the sparsity structure of the matrix.
As a result, depending on the sparsity structure of the matrix, the algorithms may not provide enough parallelism to utilize a GPU. %and may run fundamentally sequential. 

% ----------------------------------------------------
\begin{figure}[t]
  \centerline{
  \begin{subfigure}[b]{\linewidth}
  \begin{center}
  \fbox{\begin{minipage}[h!]{.75\linewidth}
    \footnotesize
    \begin{algorithmic}[1]
    \For{$s = 1,2, \dots, s_{max}$}
    \For{ $(i,j) \in S$}
     \If{$i>j$}
       \State $\ell_{ij} := (1-\omega) \ell_{ij} + \omega(a_{ij} - \sum_{k=1}^{j-1} \ell_{ik} u_{kj})/u_{ij}$
     \Else
       \State $u_{ij} := (1-\omega)u_{ij} + \omega(a_{ij} - \sum_{k=1}^{j-1} \ell_{ik} u_{kj})$
     \EndIf
    \EndFor
    \EndFor
    \end{algorithmic}
  \end{minipage}}
  \caption{FastILU to generate the sparse triangular factors $L$ and $U$ for an input matrix $A$ with a fixed sparsity structure $S$ and damping factor $\omega$.}\label{algo:fastILU_comp}
  \end{center}
  \end{subfigure}
  }
   %%
   %\quad
   %%
  \bigskip
  \centerline{
  \begin{subfigure}[b]{\linewidth}
  \begin{center}
  \fbox{\begin{minipage}[h!]{.75\linewidth}
    \footnotesize
    \begin{algorithmic}[1]
    \For{$s = 1,2,\dots, s_{max}$}
      \For{$i = 1,2,\dots,n$}
      \State $x_i^{(k)} := (b_{i} - \sum_{j=1}^{i-1} \ell_{ij} x_j^{(k-1)})/\ell_{ii}$
      \EndFor
    \EndFor
    \end{algorithmic}
  \end{minipage}}
  \caption{FastSpTRSV to compute a solution vector $x$ for a given sparse lower triangular matrix $L$ and right-hand-side vector $b$.}
  \end{center}
  \end{subfigure}
  }
  \caption{Pseudocode of Fast ILU(k) and Sparse-triangular solver with the number of sweeps given by $s_{max}$.}\label{algo:fastILU}
\end{figure}
% ----------------------------------------------------

To expose more parallelism, the
FastILU subpackage implements an iterative variant of level-based ILU factorization and sparse-triangular solver, referred to as FastILU and FastSpTRSV, respectively. The implementation utilizes Kokkos and provides portable performance across different node architectures, including NVIDIA or AMD GPUs.

FastILU is based on the algorithms proposed in~\cite{Chow:2015}.
Figure~\ref{algo:fastILU} shows the pseudocode of both FastILU and FastSpTRSV. Although each sweep of the algorithm requires about the same number of floating-point operations as the standard ILU or SpTRSV, each entry of the ILU factors or of the solution vector can be computed in parallel. Moreover, the amount of parallelism is independent of the sparsity structure of the matrix; hence, matrix reordering, which may degrade the quality of the preconditioner, is not needed to enhance the parallelism. 
Consequently, FastILU has the potential to generate an effective preconditioner efficiently without matrix ordering. Overall, if the algorithm requires a small number of sweeps to generate a preconditioner or triangular solution of a desired accuracy, then on a GPU (which can accommodate the high level of parallelism), FastILU and FastSpTRSV can be more efficient than the standard algorithms.

% ========================================================= %
% hops : srun -N1 -n1 -c1 --ntasks-per-node=4 ./Ifpack2_tif_belos.exe --xml_file=elasticity_filu.xml  --with_stacked_timer
% /nscratch/iyamaza/hops/trilinos/build/packages/ifpack2/test/belos
\begin{table}[t]\small
\centerline{
\subfloat[Standard ILU($k$).\label{tab:ilu-1}]{
\begin{tabular}{l|rrrrr}
level  & 0               & 1 & 2 & 3 & 4\\
%\# DoFs   & $2.2\times10^6$ & $8.8\times10^6$\\
\hline\hline
nnz/$n$      & 71.3 & 155.0 & 263.5 & 390.1 & 528.3\\
\hline
Numeric  & 0.019   & 0.065   & 0.225   & 0.592 & 1.381\\
Solve    & 0.052   & 0.057   & 0.061   & 0.072 & 0.115\\
\# iters & 18      & 12      & 8       & 7     & 6
\end{tabular}}
}
%
%\quad\quad\quad
%
\bigskip
\centerline{
\subfloat[FastILU($k$).]{
\begin{tabular}{l|rrrrr}
level, $k$    & 0 & 1 & 2 & 3 & 4\\
\hline\hline
\multicolumn{2}{l}{No warmup}\\
\;\;Numeric  & 0.012   & 0.045   & 0.072   & 0.124 & 0.206\\
\;\;Solve    & 0.040   & 0.063   & 0.081   & 0.099 & 0.128\\
\;\;\# iters & 18      & 14      & 13      & 12    & 12\\
\hline
\multicolumn{2}{l}{With warmup}\\
\;\;Numeric  & 0.013   & 0.047   & 0.099   & 0.212 & 0.398\\
\;\;Solve    & 0.045   & 0.059   & 0.065   & 0.070 & 0.080\\
\;\;\# iters & 18      & 13      & 9       & 7     & 6
\end{tabular}}
}
%
%\quad\quad\quad
%
\bigskip
\centerline{
\subfloat[Standard ILU($k$) with METIS.]{
\begin{tabular}{l|rrrrr}
level  & 0               & 1 & 2 & 3 & 4\\
%\# DoFs   & $2.2\times10^6$ & $8.8\times10^6$\\
\hline\hline
nnz/$n$      & 71.3 & 186.5 & 309.7 & 429.7 & 520.3\\
\hline
Numeric  & 0.017   & 0.016   & 0.054  & 1.266 & 1.316\\
Solve    & 0.055   & 0.067   & 0.077  & 0.076 & 0.081\\
\# iters & 30      & 18      & 14     & 10    & 6
\end{tabular}}
}
%
%\quad\quad\quad
%
\bigskip
\centerline{
\subfloat[FastILU($k$) with METIS.]{
\begin{tabular}{l|rrrrr}
level, $k$    & 0 & 1 & 2 & 3 & 4\\
\hline\hline
\multicolumn{2}{l}{No warmup}\\
\;\; Numeric  & 0.029   & 0.052   & 0.090  & 0.155  & 0.215\\
\;\; Solve    & 0.060   & 0.069   & 0.078  & 0.088  & 0.097\\
\;\; \# iters & 30      & 19      & 15     & 13     & 12\\
\hline
\multicolumn{2}{l}{With warmup}\\
\;\; Numeric  & 0.029   & 0.065   & 0.134  & 0.275  & 0.478\\
\;\; Solve    & 0.056   & 0.065   & 0.076  & 0.080  & 0.075\\
\;\; \# iters & 30      & 18      & 14     & 11     & 8
\end{tabular}}
}
\caption{Performance of ILU and FastILU on an NVIDIA H100 GPU for a 3D Elasticity problem on a $16^3$ grid ($n=12,288$).\label{tab:fastilu_nx16}}
\end{table}

\begin{table}[t]\small
\centerline{
\subfloat[Standard ILU($k$).\label{tab:ilu-2}]{
\begin{tabular}{l|rrrrr}
level  & 0               & 1 & 2 & 3 & 4\\
%\# DoFs   & $2.2\times10^6$ & $8.8\times10^6$\\
\hline\hline
nnz/$n$      & 76.0 & 171.4 & 302.4 & 465.2 & 655.8\\
\hline
Numeric  & 0.039   & 0.141   & 0.312   & 1.388 & 3.279\\
Solve    & 0.125   & 0.140   & 0.186   & 0.161 & 0.190\\
\# iters & 16      & 20      & 17      & 9     & 8
\end{tabular}}
}
%
%\quad\quad\quad
%
\bigskip
\centerline{
\subfloat[FastILU($k$).]{
\begin{tabular}{l|rrrrr}
level, $k$    & 0 & 1 & 2 & 3 & 4\\
\hline\hline
\multicolumn{2}{l}{No warmup}\\
\;\;Numeric  & 0.142  & 0.195 & 0.305  & 0.525 & 0.926\\
\;\;Solve    & 0.121  & 0.136 & 0.186  & 0.253 & 0.346\\
\;\;\# iters & 26     & 18    & 17     & 17    & 17\\
\hline
\multicolumn{2}{l}{With warmup}\\
\;\;Numeric  & 0.141  & 0.244 & 0.440  & 0.874 & 1.752\\
\;\;Solve    & 0.120  & 0.136 & 0.153  & 0.170 & 0.191\\
\;\;\# iters & 26     & 18    & 13     & 10    & 8
\end{tabular}}
}
%
%\quad\quad\quad
%
\bigskip
\centerline{
\subfloat[Standard ILU($k$) with METIS.]{
\begin{tabular}{l|rrrrr}
level  & 0               & 1 & 2 & 3 & 4\\
%\# DoFs   & $2.2\times10^6$ & $8.8\times10^6$\\
\hline\hline
nnz/$n$   & 76.0    & 211.5   & 364.4  & 544.1 & 676.4\\
\hline
Numeric   & 0.043   & 0.400   & 1.731  & 7.138 & 15.582\\
Solve     & 1.492   & 0.157   & 0.207  & 0.245 & 0.281\\
\# iters  & 40      & 20      & 23     & 17    & 14
\end{tabular}}
}
%
%\quad\quad\quad
%
\bigskip
\centerline{
\subfloat[FastILU($k$) with METIS.]{
\begin{tabular}{l|rrrrr}
level, $k$    & 0 & 1 & 2 & 3 & 4\\
\hline\hline
\multicolumn{2}{l}{No warmup} &\\
\;\;Numeric   & 0.137  & 0.216  & 0.384  & 0.764 & 1.189\\
\;\;Solve     & 0.103  & 0.162  & 0.226  & 0.279 & 0.346\\
\;\;\# iters  & 40     & 29     & 25     & 20    & 18\\
\hline
\multicolumn{2}{l}{With warmup}\\
\;\;Numeric  & 0.137   & 0.254  & 0.544  & 1.700 & 2.416\\
\;\;Solve    & 0.103   & 0.160  & 0.206  & 0.244 & 0.291\\
\;\;\# iters & 40      & 29     & 23     & 17    & 14
\end{tabular}}
}
\caption{Performance of ILU and FastILU on an NVIDIA H100 GPU for a 3D Elasticity problem on a $32^3$ grid ($n=98,304$). \label{tab:fastilu_nx32}}
\end{table}

Since each entry of the ILU factors is updated in parallel (in place), some of the entries used to update the entry may already have been updated by other threads (Lines 4 and 6 of Figure~\ref{algo:fastILU_comp}). Hence, the preconditioner produced by FastILU is non-deterministic, though it usually leads to a preconditioner of similar quality. 
On the other hand, FastSpTRSV updates the solution vector in a separate workspace (out of place) and calls SpMV of Kokkos-Kernels to update the solution vector at each sweep.
We did not observe variations in the number of iterations in our experiments, but the performance may be stabilized by tuning the above parameters (such as block size or $\omega$).

In addition to the standard parameters such as the level of fill,  $k$, and the number of sweeps, $s_{max}$, FastILU provides additional solver options:
\begin{itemize}
%\item Initial guess:
\item {\bf Warm up} recursively calls FastILU($k-1$) to initialize the nonzero entries in the triangular factors. This increases the cost of generating the preconditioner, but may lead to a preconditioner of higher quality.
\item {\bf Block Size} (or {\bf number of nonzeroes per thread)} allows each thread to update multiple nonzero entries in the preconditioner. It reduces the amount of parallelism, but depending on the ordering and grouping of the nonzero entries, the quality of the resulting preconditioner may improve.
\item {\bf Shift} applies the Manteuffel shift~\cite{Manteuffel:1980} before computing the preconditioner to enhance the stability of factorization.
\item {\bf Damping factor} is used to update the factor or the solution at each sweep (i.e., $\omega$ in Figure~\ref{algo:fastILU}).
\end{itemize}
The preconditioner has been used in applications such as thermal simulation~\cite{aria}. Here we show the performance results with structured 3D problems. 

Tables~\ref{tab:fastilu_nx16} and \ref{tab:fastilu_nx32} compare the performance of FastILU with the standard level-set ILU of Kokkos-Kernels~\cite{KK} on a single NVIDIA H100 GPU. The resulting ILU factors are used to precondition GMRES. We restarted the GMRES iterations after every 60 iterations, and
we considered GMRES to have converged when the relative residual
norm was reduced by six orders of magnitude.
The test matrices are for the 3D Elasticity problem on a 27-point stencil, and
the right-hand-side vector is generated such that the solution is a vector of random variables.
We used two sweeps to generate the FastILU preconditioner (i.e. $n_{max}=2$), while the standard algorithm is used to perform SpTRSV with the triangular factors generated by both Fast and standard ILUs.

For these structured 3D problems, the nested dissection reordering was not effective: with matrix reordering, the number of fills in the computed ILU factors, and the factorization time, increased, while the quality of the preconditioner degraded,
increasing the number of GMRES iterations needed for the solution convergence\footnote{For other problems, the matrix ordering may significantly reduce the factorization time, without a significant increase in the required iteration count~\cite{aria}}.
Nevertheless, the performance of FastILU was less affected by the matrix ordering (even with the increased number of fills, the factorization time of FastILU did not increase as much as it did with the standard algorithm because FastILU can more efficiently exploit the parallelism and utilize the GPU's capability).

The quality of the FastILU factorization can be improved by warming up the FastILU factorization, with increase in the factorization time (see Tables~\ref{tab:fastilu_nx16} and \ref{tab:fastilu_nx32}). The quality of the preconditioner could also be improved by increasing the number of sweeps. Table~\ref{tab:fastilu_sweep} shows the performance of FastILU(3) with different numbers of sweeps.

Overall, FastILU obtained the performance benefit over the standard ILU, especially when a large value of level is needed and the computational cost of generating the preconditioner increases
(some applications require a large value of level for the ILU preconditioner to be effective).

\begin{table}[t]\small
\centerline{
\subfloat[Performance on an NVIDIA H100 GPU. With the standard ILU(3) implementation of Kokkos-Kernels, the numerical factorization took about 1.388 seconds, while the GMRES converged with 9 iterations and took about 0.161 seconds (see Table~\ref{tab:ilu-2}).]{
\begin{tabular}{l|rrrrrr}
\# sweeps, $s_{max}$    & 1 & 2 & 3 & 4 & 5 & 6\\
\hline\hline
\;\;Numeric  & 0.470 & 0.518 & 0.562 & 0.604 & 0.646 & 0.695\\
\;\;Solve    & 0.350 & 0.252 & 0.218 & 0.194 & 0.183 & 0.170\\
\;\;\# iters & 25    & 17    & 14    & 12    & 11    & 10\\
\end{tabular}
}}
%\centerline{
%\subfloat[Performance on an AMD Instinct MI250X GPU. With the standard ILU(3) implementation of Kokkos-Kernels, the numerical %factorization took about ??? seconds, while the GMRES converged with ?? iterations and took about ??? seconds ]{
%\begin{tabular}{l|rrrrrr}
%\# sweeps, $s_{max}$    & 1 & 2 & 3 & 4 & 5 & 6\\
%\hline\hline
%\;\;Numeric  & \\
%\;\;Solve    & \\
%\;\;\# iters & \\
%\end{tabular}
%}}
\caption{Performance of FastILU(3) with different numbers of sweeps (no warmup without METIS) for a 3D Elasticity problem on a $32^3$ grid ($n=98,304$).\label{tab:fastilu_sweep}}
\end{table}

\section{Final Remarks}\label{sec:conclusion}

In this paper, we discussed the recent progress in two sparse direct solvers, Basker and Tacho, and an algebraic preconditioner, FastILU, in the ShyLU-node software packages. 
It is part of the Trilinos software framework and is designed to provide scalable solution of linear systems on a single compute node (either on multicore CPUs or on a single GPU).
Within this Trilinos software stack,
ShyLU-Basker can be used as a stand-alone solver or preconditioner for global or local problems.

Since the original publications of the solvers and preconditioner, the functionalities, robustness, and performance of ShyLU-Basker have been greatly enhanced, and these solvers are actively used for real-world applications.
To demonstrate the performance of the solvers,
we have presented the performance of these sparse direct solvers for two such real-world application problems, namely, 
Basker for solving the global problems for the Xyce Circuit Simulations on multi-core Intel CPUs
and Tacho as the local subdomain solver 
in the multi-level domain decomposition preconditioner for Albany Land-Ice Simulation of Antarctica on the Perlmutter supercomputer with NVIDIA A100 GPUs at NERSC. 
We have also demonstrated the performance of FastILU as a preconditioner for solving 3D model problems on an NVIDIA H100 GPU.

The most recent stable versions of the solvers are publicly available on the Trilinos github repository~\cite{Trilinos-website} as
they continue to evolve to address the needs of efficiently solving linear systems for new real-world applications on the emerging computers each year.

\section*{Acknowledgment}
This work was supported in part by 
the U.S. Department of Energy, Office of Science, Office of Advanced Scientific Computing Research, Scientific Discovery through Advanced Computing (SciDAC) Program through the FASTMath Institute under Contract No. DE-AC02-05CH11231 and Advanced Simulation and Computing (ASC) program at Sandia National Laboratories.
We thank Heidi Thornquist from Sandia National Labs for providing Xyce circuit matrices and valuable discussions, and Jerry Watkins and Mauro Perego also from Sandia National Labs for the help with the Albany simulations.
We also like to acknowledge the original developers of ShyLU-node including Joshua Dennis Booth, Kyungjoo Kim, and Aftab Patel.

Sandia National Laboratories is a multimission laboratory managed and operated by National Technology and Engineering Solutions of Sandia, LLC, a wholly owned subsidiary of Honeywell International, Inc., for the U.S. Department of Energy's National Nuclear Security Administration under contract DE-NA-0003525. This paper describes objective technical results and analysis. Any subjective views or opinions that might be expressed in the paper do not necessarily represent the views of the U.S. Department of Energy or the United States Government.

\bibliographystyle{IEEEtran}
\bibliography{ref}

\end{document}